\documentclass[12pt]{article}
\usepackage{amssymb}
\usepackage{amsthm}
\usepackage{amsmath}
\usepackage{graphicx}
 \newtheorem{thm}{Theorem}[section]
 
 \newtheorem{lem}[thm]{Lemma}
 
 \theoremstyle{definition}

 \numberwithin{equation}{section}

\begin{document}
\title{BKM's Criterion and Global Weak Solutions
for Magnetohydrodynamics with Zero Viscosity}
\author{Zhen Lei\footnote{School of Mathematical Sciences, Fudan
  University, Shanghai 200433, China. {\it Email:
  leizhn@yahoo.com}}\and
  Yi Zhou\footnote{School of Mathematical Sciences, Fudan
  University, Shanghai 200433, China. {\it Email: yizhou@fudan.ac.cn}}}
\date{\today}
\maketitle

\begin{abstract}
In this paper we derive a criterion for the breakdown of classical
solutions to the incompressible magnetohydrodynamic equations with
zero viscosity and positive resistivity in $\mathbb{R}^3$. This
result is analogous to the celebrated Beale-Kato-Majda's breakdown
criterion for the inviscid Eluer equations of incompressible
fluids. In $\mathbb{R}^2$ we establish global weak solutions to
the magnetohydrodynamic equations with zero viscosity and positive
resistivity for initial data in Sobolev space $H^1(\mathbb{R}^2)$.
\end{abstract}
\textbf{Keyword:} Beale-Kato-Majda's criterion, weak solutions,
magnetohydrodynamics, zero viscosity.

\section{Introduction}

The incompressible magnetohydrodynamic equations in
$\mathbb{R}^n$, $n = 2$, 3, take the form
\begin{equation}\label{a1}
\begin{cases}
u_t + u\cdot\nabla u + \nabla\big(p + \frac{1}{2}|h|^2\big)
  = \mu\Delta u + (h \cdot \nabla)h,\\[-4mm]\\
h_t - \nabla \times (u \times h) = \nu\Delta h,\\[-4mm]\\
\nabla\cdot u = 0,\quad \nabla\cdot h = 0,
\end{cases}
\end{equation}
where $u = (u_1, \cdots, u_n)^T$ is the velocity of the flows, $h
= (h_1, \cdots, h_n)^T$ is the magnetic field, $p$ is the scalar
pressure, $\mu$ is the viscosity of the fluid which is the inverse
of the Reynolds number and $\nu$ is the resistivity constant which
is inversely proportional to the electrical conductivity constant.
The system \eqref{a1} describes the macroscopic behavior of
electrically conducting incompressible fluids (see \cite{L}).

In the extremely high electrical conductivity cases which occur
frequently in the cosmical and geophysical problems, we ignore the
resistivity to have the following partially viscous
magnetohydrodynamic system (see \cite{Chandrasekhar}):
\begin{equation}\label{a2}
\begin{cases}
u_t + u\cdot\nabla u + \nabla\big(p + \frac{1}{2}|h|^2\big)
  = \mu\Delta u + (h \cdot \nabla)h,\\[-4mm]\\
h_t - \nabla \times (u \times h) = 0,\\[-4mm]\\
\nabla\cdot u = 0,\quad \nabla\cdot h = 0.
\end{cases}
\end{equation}
In the turbulent flow regime which occurs when the Reynolds
numbers is very big, we ignore the viscosity of fluids to have the
following partially viscous magnetohydrodynamic system:
\begin{equation}\label{a3}
\begin{cases}
u_t + u\cdot\nabla u + \nabla\big(p + \frac{1}{2}|h|^2\big)
  = (h \cdot \nabla)h,\\[-4mm]\\
h_t - \nabla \times (u \times h) = \nu\Delta h,\\[-4mm]\\
\nabla\cdot u = 0,\quad \nabla\cdot h = 0.
\end{cases}
\end{equation}
The local well-posedness of the Cauchy problem of the partially
viscous magnetohydrodynamic systems \eqref{a2} and \eqref{a3} is
rather standard and similar to the case of fully viscous
magnetohydrodynamic system which is done in \cite{SermangeTemam}.
At present, there is no global-in-time existence theory for strong
solutions to systems \eqref{a2} and \eqref{a3}. The notable
difference between \eqref{a2} or \eqref{a3} and its Newtonian
counterpart, the incompressible Euler equations, is that for
\eqref{a2} or \eqref{a3}, global-in-time existence has not been
established even in two dimensions (global-in-time existence is
only established in the fully viscous case $\mu > 0$, $\nu > 0$ in
two dimensions, see \cite{SermangeTemam}), even with small initial
data. The question of spontaneous apparition of singularity from a
local classical solution is a challenging open problem in the
mathematical fluid mechanics, which is similar as the cases of
ideal magnetohydrodynamics and fully viscous magnetohydrodynamics.
We just refer some of the studies on the finite time blow-up
problem in the ideal magnetohydrodynamics (see \cite{CKS, Grauer,
CordobaMarliani, CordobaFefferman} and references therein).

In the absence of a well-posedness theory, the development of
blowup/non-blowup theory is of major importance for both
theoretical and practical purposes. For incompressible Euler and
Navier-Stokes equations, the well-known Beale-Kato-Majda's
criterion \cite{BKM} says that any solution $u$ is smooth up to
time $T$ under the assumption that $\int_0^T\|\nabla\times u(t,
\cdot)\|_{L^\infty}dt < \infty$. Beale-Kato-Majda's criterion is
slightly improved by Kozono-Taniuchi \cite{KT} under the
assumption \eqref{a5}. Recently, a logrithmically improved
Beale-Kato-Majda's criterion is proven by Zhou and Lei
\cite{ZhouLei}. Caflisch-Klapper-Steele \cite{CKS} extended the
Beale-Kato-Majda's criterion to the 3D ideal magnetohydrodynamic
equations, under the assumption on both velocity field and
magnetic field: $\int_0^T\big(\|\nabla\times u(t,
\cdot)\|_{L^\infty} + \|\nabla\times h(t,
\cdot)\|_{L^\infty}\big)dt < \infty$. Motivated by numerical
experiments \cite{Hasegawa, PPS} which seem to indicate that the
velocity field plays the more important role than the magnetic
field in the regularity theory of solutions to the
magnetohydrodynamic equations, a lot of work are focused on the
regularity problem of magnetohydrodynamic equations under
assumptions only on velocity field, but not on magnetic field
\cite{HeXin1, HeXin2, Miao}. Especially, for fully viscous
magnetohydrodynamic equations, the analogy of Beale-Kato-Majda's
criterion is studied by Chen-Miao-Zhang \cite{Miao} under the
assumption only on the vorticity of velocity field, where the
authors made clever use of Littlewood-Paley theory to avoid using
the logarithmic Sobolev inequality, which successfully avoids
assuming any control of the magnetic field.

In this paper, we establish the analogous Beale-Kato-Majda's
criterion to magnetohydrodynamic equations \eqref{a3} even in the
case of zero viscosity. The partially viscous magnetohydrodynamic
system \eqref{a2} will be studied with an incompressible
viscoelastic fluid system of the Oldroyd type in our forthcoming
paper \cite{LeiMaZhou}.  The first goal of this paper is to prove
that

\begin{thm}\label{thm2}
Let $T > 0$ and $u_0, h_0 \in H^s(\mathbb{R}^n)$ for $s \geq 3$
and $n = 2$ or 3. Suppose that $(u, h)$ is a smooth solution to
the partially viscous magnetohydrodynamic system \eqref{a3} with
initial data $u(0, x) = u_0$, $h(0, x) = h_0$. Then $(u, h)$ is
smooth up to time $T$ provided that
\begin{equation}\label{a5}
\int_0^T\|\nabla\times u(t, \cdot)\|_{{\rm BMO}}dt < \infty.
\end{equation}
\end{thm}

Our key observation is that under the assumption \eqref{a5}, if we
start from the time $T_\star < T$ which is as close to $T$ as
possible, we find that
\begin{equation}\nonumber
\|(u, h)(t, \cdot)\|_{H^1} \leq C_\star\sup_{T_\star \leq s \leq
t}\|(u, h)\|_{H^3}^\delta,\quad \delta\ {\rm is\ as\ small\ as\
we\ want}, T_\star \leq t < T.
\end{equation}
From this point of view, the nonlinearity of velocity field $u$
and magnetic field $h$ can be as weak as what we want in the sense
of $H^1$ norm. By using this observation, we find that logarithmic
Sobolev inequalities can still be used to prove Theorem
\ref{thm2}.

Our second result is the global existence and regularity of weak
solutions to partially viscous magnetohydrodynamic system
\eqref{a3} in $\mathbb{R}^2$:
\begin{thm}\label{thm3}
Let $u_0, h_0 \in H^1(\mathbb{R}^2)$. Then there exists a global
weak solution $(u, h)$ with $u \in L^\infty(0, \infty; H^1)$ and
$h \in L^\infty(0, \infty; H^1) \cap L^2(0, \infty; H^{2})$ to the
partially viscous magnetohydrodynamic system \eqref{a3} with
initial data $u(0, x) = u_0$, $h(0, x) = h_0$.
\end{thm}

The observation of the above result lies in that we can obtain the
$H^1$ estimate of the velocity field and magnetic field. Then the
proof of Theorem \ref{thm3} is elementary and standard. Thus, we
will just present the $H^1$ estimate of the velocity field and
magnetic field in this paper, leaving the establishment of global
weak solutions to interesting readers (for example, see the books
\cite{Temam, Majda} as references). See section 4 for more
details.

The rest of this paper is organized as follows: In section 2 we
recall some inequalities which will be used in this paper. In
section 3 we derive the Beale-Kato-Majda's criterion for the
partially viscous magnetohydrodynamic equaitons \eqref{a3} and
present the proof of Theorem \ref{thm2}. Then in section 4 we
prove the $H^1$ \textit{a priori} estimate of the velocity field
$u$ and the magnetic field $h$ in $\mathbb{R}^2$ which implies the
result in Theorem \ref{thm3}.

\section{Preliminaries}

First of all, let us recall the following multiplicative
inequalities.
\begin{lem}\label{lem1}
The following interpolation inequalities holds.
\begin{equation}\label{b1}
\begin{cases}
{\rm in\ \mathbb{R}^2}:
\begin{cases}
\|f\|_{L^\infty} \leq C_0\|f\|_{L^2}^{\frac{1}{2}}\|\nabla^2
  f\|_{L^2}^{\frac{1}{2}},\\[-4mm]\\
\|f\|_{L^4} \leq C_0\|f\|_{L^2}^{\frac{1}{2}}\|\nabla
  f\|_{L^2}^{\frac{1}{2}},\\[-4mm]\\
\|f\|_{L^4} \leq C_0\|f\|_{L^2}^{\frac{3}{4}}
  \|\nabla^2 f\|_{L^2}^{\frac{1}{4}},\\[-4mm]\\
\end{cases}\\
{\rm in\ \mathbb{R}^3}:
\begin{cases}
\|f\|_{L^\infty} \leq C_0\|f\|_{L^2}^{\frac{1}{4}}\|\nabla^2
  f\|_{L^2}^{\frac{3}{4}},\\[-4mm]\\
\|f\|_{L^4} \leq C_0\|f\|_{L^2}^{\frac{5}{8}}
  \|\nabla^2 f\|_{L^2}^{\frac{3}{8}},\\[-4mm]\\
\|\nabla f\|_{L^2} \leq C_0\|f\|_{L^2}^{\frac{1}{2}}\|\nabla^2
  f\|_{L^2}^{\frac{1}{2}},
\end{cases}
\end{cases}
\end{equation}
where $C_0$ is an absolute positive constant.
\end{lem}
\begin{proof}
The above inequalities are of course well-known. In fact, they can
just proved by Sobolev embedding theorems and the scaling
techniques. As an example, let us present the proof of the first
inequalities. By Sobolev embedding theorem $H^2(\mathbb{R}^2)
\hookrightarrow L^\infty(\mathbb{R}^2)$, one has
\begin{equation}\nonumber
\|f\|_{L^\infty} \leq C_0\big(\|f\|_{L^2} + \|\nabla^2
f\|_{L^2}\big).
\end{equation}
The constant $C_0$ in above inequality is independent of $f \in
H^2(\mathbb{R}^2)$. Thus, for any given $0 \neq f \in
H^2(\mathbb{R}^2)$ and  any $\lambda > 0$, let us define
$f^\lambda(x) = f(\lambda x)$. Then one  has
\begin{equation}\nonumber
\|f^\lambda\|_{L^\infty} \leq C_0\big(\|f^\lambda\|_{L^2} +
\|\nabla^2 f^\lambda\|_{L^2}\big),
\end{equation}
which is equivalent to
\begin{equation}\nonumber
\|f\|_{L^\infty} \leq C_0\big(\lambda^{-1}\|f\|_{L^2} + \lambda
\|\nabla^2 f\|_{L^2}\big).
\end{equation}
Taking $\lambda = (\|f\|_{L^2}\|\nabla^2
f\|_{L^2})^{\frac{1}{2}}$, one gets the first inequality in
\eqref{b1}.
\end{proof}

Next, let us recall the following well-known inequalities. In
fact, the first one is Gagliardo-Nirenberg inequality and the
second one is a direct consequence of the chain rules and
Gagliardo-Nirenberg inequality.
\begin{lem}\label{lem2}
The following inequalities holds:
\begin{equation}\label{b2}
\begin{cases}
\|\nabla^iu\|_{L^{\frac{2s}{i}}} \leq C_0\|u\|_{L^\infty}^{1 -
  \frac{i}{s}}\|\nabla^su\|_{L^2}^{\frac{i}{s}},
  \quad 0 \leq i \leq s,\\[-4mm]\\
\big\|\nabla^s(u\cdot\nabla u) - u\cdot\nabla\nabla^s u\big\| \leq
  C_0\|\nabla u\|_{L^\infty}\|\nabla^s u\|_{L^2},\ s \geq 1.
\end{cases}
\end{equation}
\end{lem}

At last, let us recall the following logarithmic Sobolev
inequality which is proved in \cite{KT} and is an improved version
of that in \cite{BKM} (see also \cite{Brezis, Brezis2}).
\begin{lem}\label{lem4}
Let $n = 2$ or 3 and $p > n$. The following logarithmic Sobolev
embedding theorem holds for all divergence free vector fields:
\begin{equation}\label{b4}
\|\nabla f\|_{L^\infty(\mathbb{R}^n)} \leq C_0\big[1 +
\|f\|_{L^2(\mathbb{R}^n)} + \|\nabla\times f\|_{{\rm
BMO}(\mathbb{R}^n)}\ln\big(1 + \|f\|_{W^{2,
p}(\mathbb{R}^n)}\big)\big].
\end{equation}
\end{lem}

\section{Beale-Kato-Majda's Criterion for Magnetohydrodynamic
Equations with Zero Viscosity}

In this section we prove the analogous Beale-Kato-Majda's
criterion for the partially viscous magnetohydrodynamic equations
\eqref{a3} in $\mathbb{R}^2$ and $\mathbb{R}^3$. First of all, for
classical solutions to \eqref{a3}, one has the following basic
energy law
\begin{eqnarray}\label{c1}
\frac{1}{2}\frac{d}{dt}\big(\|u\|_{L^2}^2 + \|h\|_{L^2}^2\big) +
\nu\|\nabla h\|_{L^2}^2 = 0.
\end{eqnarray}

For any integer $s \geq 1$, applying $\nabla^s$  to the velocity
field equation and magnetic field equation, and taking the $L^2$
inner product of the resulting equations with $\nabla^su$ and
$\nabla^sh$ respectively, one has
\begin{eqnarray}\label{d}
&&\frac{1}{2}\frac{d}{dt}\big(\|\nabla^s u\|_{L^2}^2 +
  \|\nabla^s h\|_{L^2}^2\big) + \nu\|\nabla^s\nabla
  h\|_{L^2}^2\\\nonumber
&&= - \int_{\mathbb{R}^n}\nabla^su\big[\nabla^s(u\cdot\nabla u) -
  u\cdot\nabla\nabla^s u\big]dx\\\nonumber
&&\quad -\ \int_{\mathbb{R}^n}\nabla^sh\big[\nabla^s(u\cdot\nabla
  h) - u\cdot\nabla\nabla^sh\big]dx\\\nonumber
&&\quad +\ \int_{\mathbb{R}^n}\nabla^su\big[\nabla^s(h\cdot\nabla
  h) - h\cdot\nabla\nabla^sh\big]dx\\\nonumber
&&\quad +\ \int_{\mathbb{R}^n}\nabla^sh\big[\nabla^s(h\cdot\nabla
  u) - h\cdot\nabla\nabla^s u\big]dx,
\end{eqnarray}
where we have used the divergence free condition $\nabla\cdot u =
\nabla\cdot h = 0$.

For simplicity, we just set $\nu = 1$ below. Let $s = 1$ in
\eqref{d}. One can easily derive the following estimate:
\begin{eqnarray}\nonumber
&&\frac{1}{2}\frac{d}{dt}\big(\|\nabla u\|_{L^2}^2 +
  \|\nabla h\|_{L^2}^2\big) + \|\nabla^2 h\|_{L^2}^2\\\nonumber
&&\leq C\|\nabla u\|_{L^\infty}\big(\|\nabla u\|_{L^2}^2 +
\|\nabla h\|_{L^2}^2\big)
\end{eqnarray}
which gives that
\begin{eqnarray}\label{d1}
&&\|\nabla u(t, \cdot)\|_{L^2}^2 + \|\nabla h(t, \cdot)\|_{L^2}^2
  + \int_{t_0}^t\|\nabla^2 h(s, \cdot)\|_{L^2}^2ds\\\nonumber
&&\leq C\big(\|\nabla u(t_0, \cdot)\|_{L^2}^2 + \|\nabla
  h(t_0, \cdot)\|_{L^2}^2\big)\exp\Big\{C\int_{t_0}^t\|\nabla
  u\|_{L^\infty}ds\Big\}.
\end{eqnarray}
Noting \eqref{a5}, one concludes that for any small constant
$\epsilon > 0$, there exists $T_\star < T$ such that
\begin{equation}\label{d2}
\int_{T_\star}^T\|\nabla\times u\|_{{\rm BMO}}ds \leq \epsilon.
\end{equation}
Let us denote
\begin{eqnarray}\label{d3}
M(t) = \sup_{T_\star \leq s \leq t}\big(\|\nabla^3 u(s,
\cdot)\|_{L^2}^2 + \|\nabla^3h(s, \cdot)\|_{L^2}^2\big),\quad
T_\star \leq t < T.
\end{eqnarray}
By \eqref{c1}, \eqref{d1}, \eqref{d2}, \eqref{d3} and Lemma
\ref{lem4},  one has
\begin{eqnarray}\label{d4}
&&\|\nabla u(t, \cdot)\|_{L^2}^2 + \|\nabla h(t, \cdot)\|_{L^2}^2
  + \int_{T_\star}^t\|\nabla^2 h(s, \cdot)\|_{L^2}^2ds\\\nonumber
&&\leq C_\star\exp\Big\{C_0\int_{T_1}^t\|\nabla\times u
  \|_{{\rm MBO}}\ln\big(1 + \|u\|_{H^3}^2 +
  \|h\|_{H^3}^2\big)ds\Big\}\\\nonumber
&&\leq C_\star\exp\big\{C_0\epsilon\ln\big(1 + M(t)\big)\big\} =
  C_\star \big(1 + M(t)\big)^{C_0\epsilon},\quad T_\star \leq t < T,
\end{eqnarray}
where $C_\star$ depends on $\|\nabla u(T_\star, \cdot)\|_{L^2}^2 +
\|\nabla h(T_\star, \cdot)\|_{L^2}^2$, while $C_0$ is an absolute
positive constant given in Lemma \ref{lem4}. We remark here that
in the case of $n = 2$, one can just set $\epsilon = 0$. See
section 4 for more details.

For simplicity below, we will set $s = 3$ and present the estimate
of the right hand side of \eqref{d}. First of all, by Lemma
\ref{lem2}, it is easy to see that
\begin{eqnarray}\label{d5}
\Big|\int_{\mathbb{R}^n}\nabla^3u\big[\nabla^3(u\cdot\nabla u) -
u\cdot\nabla\nabla^3 u\big]dx\Big| \leq C\|\nabla
u\|_{L^\infty}\|\nabla^3 u\|_{L^2}^2.
\end{eqnarray}
Next, by integrating by parts, one has
\begin{eqnarray}\nonumber
&&\Big|\int_{\mathbb{R}^n}\nabla^3h\big[\nabla^3(u\cdot\nabla
  h) - u\cdot\nabla\nabla^3 h\big]dx\Big|\\\nonumber
&&\quad +\ \Big|\int_{\mathbb{R}^n}\nabla^3h
  \big[\nabla^3(h\cdot\nabla u) - h\cdot\nabla
  \nabla^3 u\big]dx\Big|\\\nonumber
&&\leq 4\|\nabla u\|_{L^\infty} \|\nabla^3h\|_{L^2}^2\\\nonumber
&&\quad +\ 3\Big|\int_{\mathbb{R}^n} \nabla^3h
  \big(\nabla^2u\cdot\nabla\nabla h\big)dx\Big|
  + \Big|\int_{\mathbb{R}^n}\nabla^3h\big(\nabla^3
  u\cdot\nabla h\big)dx\Big|\\\nonumber
&&\quad +\ 3\Big|\int_{\mathbb{R}^n}\nabla^3h\nabla^2
  h\cdot\nabla\nabla udx\Big| + 3\Big|\int_{\mathbb{R}^n}\nabla^3h\nabla
  h\cdot\nabla\nabla^2 udx\Big|\\\nonumber
&&\leq 14\|\nabla u\|_{L^\infty} \|\nabla^3h\|_{L^2}^2
  + 10\|\nabla u\|_{L^\infty}\|\nabla^2h\|_{L^2}
  \|\nabla^4h\|_{L^2}\\\nonumber
&&\quad +\ 4\|\nabla^2 u\|_{L^4}\|\nabla
  h\|_{L^4}\|\nabla^4h\|_{L^2}.
\end{eqnarray}
By Lemma \ref{lem1} and Lemma \ref{lem2}, we estimate
\begin{eqnarray}\nonumber
&&10\|\nabla u\|_{L^\infty}\|\nabla^2h\|_{L^2}
  \|\nabla^4h\|_{L^2}\\\nonumber
&&\leq \frac{1}{8}\|\nabla^4h\|_{L^2}^2 + C\|\nabla u
  \|_{L^\infty}^2\|\nabla^2 h\|_{L^2}^2\\\nonumber
&&\leq \frac{1}{8}\|\nabla^4h\|_{L^2}^2 + C\|\nabla u
  \|_{L^\infty}\|\nabla u\|_{L^2}^{\frac{1}{4}}\|\nabla^3
  u\|_{L^2}^{\frac{3}{4}}\|\nabla h\|_{L^2}\|\nabla^3 h\|_{L^2}\\\nonumber
&&\leq \frac{1}{8}\|\nabla^4h\|_{L^2}^2 + C_\star\|\nabla u
  \|_{L^\infty}M(t)^{\frac{7}{8}}\big(1 +
  M(t)\big)^{\frac{5C_0\epsilon}{8}}
\end{eqnarray}
in 3D case and
\begin{eqnarray}\nonumber
&&10\|\nabla u\|_{L^\infty}\|\nabla^2h\|_{L^2}
  \|\nabla^4h\|_{L^2}\\\nonumber
&&\leq \frac{1}{8}\|\nabla^4h\|_{L^2}^2 + C\|\nabla u
  \|_{L^\infty}\|\nabla u\|_{L^2}^{\frac{1}{2}}\|\nabla^3
  u\|_{L^2}^{\frac{1}{2}}\|\nabla h\|_{L^2}\|\nabla^3 h\|_{L^2}\\\nonumber
&&\leq \frac{1}{8}\|\nabla^4h\|_{L^2}^2 + C_\star\|\nabla u
  \|_{L^\infty}M(t)^{\frac{3}{4}}\big(1 +
  M(t)\big)^{\frac{3C_0\epsilon}{4}}
\end{eqnarray}
in 2D case, where we used \eqref{d4}. On the other hand, one can
similarly estimate
\begin{eqnarray}\nonumber
&&4\|\nabla^2 u\|_{L^4}\|\nabla h\|_{L^4}
  \|\nabla^4h\|_{L^2}\\\nonumber
&&\leq 4\|\nabla u\|_{L^\infty}^{\frac{1}{2}}\|\nabla^3 u
  \|_{L^2}^{\frac{1}{2}}\|\nabla h\|_{L^2}^{\frac{5}{8}}
  \|\nabla^3 h\|_{L^2}^{\frac{3}{8}}\|\nabla^4h\|_{L^2}\\\nonumber
&&\leq \frac{1}{8}\|\nabla^4h\|_{L^2}^2 + C\|\nabla u\|_{L^\infty}
  \|\nabla h\|_{L^2}^{\frac{5}{4}}\|\nabla^3 u\|_{L^2}
  \|\nabla^3 h\|_{L^2}^{\frac{3}{4}}\\\nonumber
&&\leq \frac{1}{8}\|\nabla^4h\|_{L^2}^2 + C_\star\|\nabla
  u\|_{L^\infty}M(t)^{\frac{7}{8}}\big(1
  + M(t)\big)^{\frac{5C_0\epsilon}{8}}
\end{eqnarray}
in 3D case and
\begin{eqnarray}\nonumber
&&4\|\nabla^2 u\|_{L^4}\|\nabla h\|_{L^4}
  \|\nabla^4h\|_{L^2}\\\nonumber
&&\leq 4\|\nabla u\|_{L^\infty}^{\frac{1}{2}}\|\nabla^3 u
  \|_{L^2}^{\frac{1}{2}}\|\nabla h\|_{L^2}^{\frac{3}{4}}
  \|\nabla^3 h\|_{L^2}^{\frac{1}{4}}\|\nabla^4h\|_{L^2}\\\nonumber
&&\leq \frac{1}{8}\|\nabla^4h\|_{L^2}^2 + C_\star\|\nabla
  u\|_{L^\infty}M(t)^{\frac{3}{4}}\big(1
  + M(t)\big)^{\frac{3C_0\epsilon}{4}}
\end{eqnarray}
in 2D case. Consequently, one has
\begin{eqnarray}\label{d6}
&&\Big|\int_{\mathbb{R}^n}\nabla^sh\big[\nabla^s(u\cdot\nabla
  h) - u\cdot\nabla\nabla^s h\big]dx\Big|\\\nonumber
&&\quad +\ \Big|\int_{\mathbb{R}^n}\nabla^sh
  \big[\nabla^s(h\cdot\nabla u) - h\cdot\nabla
  \nabla^s u\big]dx\Big|\\\nonumber
&&\leq \frac{1}{4}\|\nabla^4h\|_{L^2}^2 + C_\star\|\nabla
  u\|_{L^\infty}\big[1 + M(t)\big]
\end{eqnarray}
provided that
\begin{eqnarray}\label{d7}
\epsilon \leq \frac{1}{5C_0}.
\end{eqnarray}

It remains to estimate the last term on the right hand side of
\eqref{d}. By integrating by parts:
\begin{eqnarray}\nonumber
\Big|\int_{\mathbb{R}^n}\nabla^3u\nabla^2h\cdot\nabla
  \nabla hdx\Big| \leq \Big|\int_{\mathbb{R}^n}\nabla^2u\nabla^3h\cdot\nabla
  \nabla hdx\Big| + \Big|\int_{\mathbb{R}^n}\nabla^2u\nabla^2h\cdot\nabla
  \nabla^2 hdx\Big|,
\end{eqnarray}
one similarly has the following estimate
\begin{eqnarray}\label{d8}
&&\Big|\int_{\mathbb{R}^n}\nabla^3u\big[\nabla^3(h\cdot\nabla
  h) - h\cdot\nabla\nabla^3 h\big]dx\Big|\\\nonumber
&&\leq \Big|\int_{\mathbb{R}^n}\nabla^3u\nabla^3h\cdot\nabla
  hdx\Big| + 3\Big|\int_{\mathbb{R}^n}\nabla^3u\nabla^2h\cdot\nabla
  \nabla hdx\Big|\\\nonumber
&&\quad +\ 3\Big|\int_{\mathbb{R}^n}\nabla^3u\nabla
  h\cdot\nabla\nabla^2hdx\Big|\\\nonumber
&&\leq \frac{1}{4}\|\nabla^4h\|_{L^2}^2 + C\|\nabla u\|_{L^\infty}
  \big[1 + M(t)\big].
\end{eqnarray}

Combining \eqref{d} with \eqref{d5}, \eqref{d6} and \eqref{d8} and
using Lemma \ref{lem4}, we arrive at
\begin{eqnarray}\nonumber
&&\frac{d}{dt}\big(\|\nabla^3 u\|_{L^2}^2 +
  \|\nabla^3 h\|_{L^2}^2\big) + \|\nabla^3\nabla
  h\|_{L^2}^2\\\nonumber
&&\leq C_\star\|\nabla u\|_{L^\infty}\big[1 +
  M(t)\big]\\\nonumber
&&\leq C_\star\|\nabla\times u\|_{{\rm BMO}}\ln\big(1 +
M(t)\big)\big[1 + M(t)\big]
\end{eqnarray}
for all $T_\star \leq t < T$. Integrating the above inequality
with respect to time from $T_\star$ to $t \in [s , T)$ and using
\eqref{d2}, we have
\begin{eqnarray}\nonumber
&&1 + \|\nabla^3 u(s, \cdot)\|_{L^2}^2 +
  \|\nabla^3 h(s, \cdot)\|_{L^2}^2\\\nonumber
&&\leq 1 + \|\nabla^3 u(T_\star, \cdot)\|_{L^2}^2 +
  \|\nabla^3 h(T_\star, \cdot)\|_{L^2}^2\\\nonumber
&&\quad +\  C_\star\int_{T_\star}^s\|\nabla\times u(\tau,
  \cdot)\|_{{\rm BMO}}\big[1 + M(\tau)\big]\ln\big(1 +
  M(\tau)\big)d\tau,
\end{eqnarray}
which implies
\begin{eqnarray}\nonumber
&&1 +  M(t) \leq 1 + \|\nabla^3 u(T_\star, \cdot)\|_{L^2}^2 +
  \|\nabla^3 h(T_\star, \cdot)\|_{L^2}^2\\\nonumber
&&\quad +\  C_\star\int_{T_\star}^t\|\nabla\times u(\tau,
  \cdot)\|_{{\rm BMO}}\big[1 + M(\tau)\big]\ln\big(1 +
  M(\tau)\big)d\tau.
\end{eqnarray}
Then gronwall's inequality gives
\begin{eqnarray}\label{d9}
&&1 + \|\nabla^3 u(t, \cdot)\|_{L^2}^2 + \|\nabla^3 h(t,
  \cdot)\|_{L^2}^2\\\nonumber
&&\leq \big(1 + \|\nabla^3 u(T_\star, \cdot)\|_{L^2}^2 +
  \|\nabla^3 h(T_\star, \cdot)\|_{L^2}^2\big)\exp\exp\{C_\star\epsilon\}
\end{eqnarray}
for all $T_\star \leq t < T$. Noting that the right hand side of
\eqref{d9} is independent of $t$ for $T_\star \leq t < T$, one
concludes that \eqref{d9} is also valid for $t = T$ which means
that $\big(u(T, \cdot), h(T, \cdot)\big) \in H^3(\mathbb{R}^n)$.

\section{Global Weak Solutions to the Magnetohydrodynamic
Equations with Zero Viscosity in 2D}

In this section we prove Theorem \ref{thm3}. To this end,
following a standard procedure, we establish an approximate system
of \eqref{a3} with smoothed initial data which admits a unique
local classical solution. We will establish the global \textit{a
priori} $H^1$-estimate for the approximate system in terms of the
$H^1$-norm of the original initial data $u_0$ and $h_0$, which
implies that the magnetohydrodynamic equations \eqref{a3} in
$\mathbb{R}^2$ possess a global weak solution. Below we will just
present the global \textit{a priori} estimate for classical
solutions to the magnetohydrodynamic equations \eqref{a3} with $n
= 2$, while we refer the reader to as \cite{Temam, Majda} for the
standard procedure to establish the global weak solutions and to
construct approximate systems.

Now let us denote ${\rm curl}(v) = \nabla \times v =
\partial_{x_1}v_2 - \partial_{x_2}v_1$ for two dimensional vector $v$.
Applying ${\rm curl}$ to \eqref{a3}, one has
\begin{equation}\label{c8}
\begin{cases}
(\nabla\times u)_t + (u\cdot\nabla)(\nabla\times u)
  = (h \cdot\nabla)(\nabla\times h),\\[-4mm]\\
(\nabla\times h)_t + (u\cdot\nabla)(\nabla\times h)
  = \nu\Delta(\nabla\times h)\\
\quad\quad   + (h\cdot\nabla)(\nabla\times u) + 2{\rm tr}(\nabla
  u\nabla^\perp h).
\end{cases}
\end{equation}
Taking the $L^2$ inner product of the above equations with
$\nabla\times u$ and $\nabla\times h$ respectively, and noting the
interpolation inequality in Lemma \ref{lem1}, one has
\begin{eqnarray}\label{c2}
&&\frac{1}{2}\frac{d}{dt}\big(\|\nabla \times u\|_{L^2}^2 +
  \|\nabla \times h\|_{L^2}^2\big) + \nu\|\nabla(\nabla \times
  h)\|_{L^2}^2\\\nonumber
&&= \int_{\mathbb{R}^2}\big[(\nabla\times u)(h \cdot\nabla)
  (\nabla\times h) + (\nabla\times h)(h\cdot\nabla)
  (\nabla\times u)\big]dx\\\nonumber
&&\quad +\ 2\int_{\mathbb{R}^2}(\nabla\times h){\rm
  tr}(\nabla u\nabla^\perp h)dx\\\nonumber
&&\leq 0 + 2\|\nabla h\|_{L^4}^2\|\nabla u\|_{L^2} \leq
  C\|\nabla h\|_{L^2}\|\nabla u\|_{L^2}\|\Delta h\|_{L^2}\\\nonumber
&&\leq \frac{\nu}{2}\|\nabla^\perp(\nabla\times h)\|_{L^2}^2 +
  \frac{C}{\nu}\|\nabla h\|_{L^2}^2\|\nabla\times u\|_{L^2}^2,
\end{eqnarray}
where we used Calderon-Zygmund theory and the fact that $\Delta h
= \nabla^\perp(\nabla\times h)$. Using the basic energy law
\eqref{c1}, we derive from \eqref{c2} that
\begin{eqnarray}\label{c3}
&&\|\nabla \times u\|_{L^2}^2 + \|\nabla \times
  h\|_{L^2}^2 + \nu\int_0^t\|\nabla (\nabla \times h)\|_{L^2}^2ds\\\nonumber
&&\leq \big(\|\nabla \times u_0\|_{L^2}^2 + \|\nabla \times
  h_0\|_{L^2}^2\big)\exp\Big\{\frac{2C}{\nu}\int_0^t\|\nabla
  h\|_{L^2}^2ds\Big\}\\\nonumber
&&\leq \big(\|\nabla \times u_0\|_{L^2}^2 + \|\nabla \times
  h_0\|_{L^2}^2\big)\exp\big\{\frac{C\big(\|u_0\|_{L^2}^2 +
  \|h_0\|_{L^2}^2\big)}{\nu^2}\big\}.
\end{eqnarray}
This combining with the basic energy law \eqref{c1} gives a global
uniform bound for the velocity field $u$ and the magnetic field
$h$ in $u \in L^\infty(0, \infty; H^1)$ and $h \in L^\infty(0,
\infty; H^1) \cap L^2(0, \infty; H^2)$.

\bigskip

\textit{Final Remark: In fact, one may improve the regularity of
the weak solutions if the initial data is more regular. For
example, for some $p, q \in (2, \infty)$, using the estimates for
linear Stokes system (see \cite{GigaSohr}) and transport
equations, one can get a uniform bound for
$\|\nabla^2h\|_{L^q_t(L^p_x)} + \|\nabla\times
u\|_{L^\infty_t(L^p_x)}$. However, at the moment we are not able
to get the global classical solutions of the partially viscous
magnetohydrodynamic equations \eqref{a3} for $n = 2$. We will
investigate this issue further in our future work.}

\bigskip

\section*{Acknowledgments}
The authors would like to thank professors Congming Li and N.
Masmoudi for constructive discussions. Zhen Lei was partially
supported by National Science Foundation of China (grant no.
10801029) and a special Postdoctoral Science Foundation of China
(grant no. 200801175). Yi Zhou is partially supported by the
National Science Foundation of China under grants 10728101, the
973 project of the Ministry of Science and Technology of China,
the Doctoral Program Foundation of the Ministry of Education of
China and the ¡±111¡± project (B08018).

\bibliographystyle{amsplain}

\end{document}